\documentclass{amsart}
\usepackage{graphicx}
\usepackage{amssymb,amstext,amsmath,amscd,amsthm,amsfonts,enumerate,graphicx,latexsym,stmaryrd,multicol,tikz-cd,hyperref,setspace,cleveref,color,comment}
\usepackage[all]{xy}
\usepackage[colorinlistoftodos]{todonotes}

\numberwithin{equation}{section}

\newtheorem{theorem}{Theorem}[section]

\newtheorem{lemma}[theorem]{Lemma}
\newtheorem{proposition}[theorem]{Proposition}
\newtheorem{corollary}[theorem]{Corollary}

\theoremstyle{definition}
\newtheorem{definition}[theorem]{Definition}
\newtheorem{remark}[theorem]{Remark}
\newtheorem{notation}[theorem]{Notation}
\newtheorem{example}[theorem]{Example}

\newcommand\Spec{\operatorname{Spec}}
\newcommand\Hom{\operatorname{Hom}}
\newcommand\Ext{\operatorname{Ext}}
\newcommand\Tor{\operatorname{Tor}}

\newcommand\depth{\operatorname{depth}}

\newcommand\coker{\operatorname{coker}}
\newcommand\pd{\operatorname{pd}}
\newcommand\id{\operatorname{id}}

\newcommand\Rfd{\operatorname{Rfd}}
\newcommand\Gid{\operatorname{Gid}}
\newcommand\Gpd{\operatorname{Gpd}}
\newcommand\Gfd{\operatorname{Gfd}}
\newcommand\gdim{\operatorname{G-dim}}

\newcommand\Cid{\operatorname{CI-id}}
\newcommand\cdim{\operatorname{CI-dim}}
\newcommand\uCid{\operatorname{CI*-id}}
\newcommand\hdim{\operatorname{H-dim}}
\newcommand\hid{\operatorname{Hid}}

\newcommand\px{\mathfrak{p}}
\newcommand\mx{\mathfrak{m}}

\newcommand{\sD}{\operatorname{\mathsf{D}}}
\newcommand{\rb}{\square}
\newcommand{\rf}{\mathsf{f}}
\newcommand{\Dbf}{\sD_\rb^\rf}
\newcommand{\Db}{\sD_\rb}
\newcommand{\Dlb}{\sD_{\sqsubset}}
\newcommand{\Dlbf}{\sD_{\sqsubset}^\rf}
\newcommand{\Drb}{\sD_{\sqsupset}}
\newcommand{\Drbf}{\sD_{\sqsupset}^\rf}
\newcommand{\ldt}{\otimes^\mathbf{L}}

\newcommand\RHom{\operatorname{\mathbf{R}Hom}}

\newcommand\HH{\operatorname{H}}
\newcommand\CC{\operatorname{C}}
\newcommand\cone{\operatorname{cone}}

\begin{document}

\title[]{Finiteness of homological dimensions of RHom complexes}

\author{Paulo Martins}
 \address{Universidade de S{\~a}o Paulo -
 ICMC, Caixa Postal 668, 13560-970, S{\~a}o Carlos-SP, Brazil}
 \email{paulomartinsmtm@gmail.com}

 \author{Victor D. Mendoza-Rubio }
 \address{Universidade de S{\~a}o Paulo -
 ICMC, Caixa Postal 668, 13560-970, S{\~a}o Carlos-SP, Brazil}
 \email{vicdamenru@gmail.com}

 \author{Zachary Nason}
 \address{Department of Mathematics, University of Nebraska, Lincoln, NE 68588-0130, USA}
 \email{znason2@huskers.unl.edu}

\keywords{derived Hom complex, complete intersection dimension, Ext modules, projective dimension, injective dimension, Gorenstein dimension, dualizing complex}
\subjclass[2020]{Primary: 13D05, 13D07, 13D09,  13H10, 13C10; Secondary: 
13D02}
\begin{abstract} This paper studies derived $\Hom$ complexes with finite projective, Gorenstein projective, or complete intersection dimension over a Noetherian local ring $R$. We establish several stability results concerning finite complete intersection homological dimension of derived $\Hom$ complexes. In addition, for $R$-complexes $M$ and $N$ we study when $\RHom_R(M,N)$ has finite projective, complete intersection dimension or Gorenstein projective dimension under either of the assumptions that $\Gpd_R(N)<\infty$ or $\Gid_R(M)<\infty$. We derive isomorphisms that mimic derived Hom evaluation and derived tensor evaluation when the derived Hom complex has finite Gorenstein projective dimension; while finite projective or complete intersection dimension of the derived Hom complex yields new stability results.

\end{abstract}

\maketitle

\section{Introduction}
Throughout this paper, unless otherwise specified, we assume that $R$ is a commutative Noetherian local ring with maximal ideal $\mx$. We will frequently regard $M$ and $N$ as complexes of (not necessarily finitely generated) $R$-modules; however, whenever we use $M$ and $N$ to denote $R$-modules, they will be assumed to be finitely generated.

In the literature, there are several results concerning finiteness of homological dimensions of $\RHom$ complexes. In particular, given $R$-complexes $M$ and $N$, there are various stability results showing that the finiteness of homological dimensions of $\RHom(M, N)$ influences the finiteness of homological dimensions of $M$ or $N$. We recall some standard stability results, which can be found in  \cite[Propositions 16.4.32, 16.4.34, 16.4.19]{DerivedCategoryMethodsInCommutativeAlgebra},  \cite[Theorem 5.6]{GdimensionOverLocalHomomorphisms} and \cite[Theorem 3.11]{CompleteIntersectionDimensionsForComplexes}. Let $M$ and $N$ be non-acyclic complexes in $\Dbf(R)$. Then the following implications hold:
\begin{enumerate}
    \item  $\id_R (\RHom_R(M,N))<\infty \Longrightarrow \pd_R( M)<\infty \text{ and }\id_R (N)<\infty$. 
    \item $\pd_R( \RHom_R(M,N)) <\infty \text{ and } \id_R (N)<\infty 
     \Longrightarrow \id_R (M)<\infty$.
    \item  $\pd_R (\RHom_R(M,N))<\infty \text{ and } \pd_R (M)<\infty 
     \Longrightarrow \pd_R (N)<\infty$.
     \item $\gdim_R (\RHom_R(M,N)) <\infty \text{ and } \pd_R (M)<\infty 
     \Longrightarrow \gdim_R (N)<\infty$.
     \item $\cdim_R (\RHom_R(M,N)) <\infty $ and $\pd_R (M)<\infty$ 
    $\Longrightarrow$ $\cdim_R (N)<\infty$.
\end{enumerate}
More stability results concerning the finiteness of homological dimensions of $\RHom$ can be found in \cite{CompleteIntersectionDimensionsandFoxbyClasses,sather2021complete,yassemi1995g,christensen2017gorenstein,Totushek}. In addition, finiteness of homological dimensions of homologies of RHom complexes are also studied; see \cite{injectivedimensiontakahashi,AuslanderReitenCinjetivadimensionVanishingofExt, NumericalApectsofComplexesOfFiniteHomologicalDimensions, FiniteHomologicalDimensionOfHomv3, OnModulesWhoseDualIsOfFiniteGorensteinDimension, HomologicalDimensionsTheGorensteinPropertyAndSpecialCasesOfSomeConjectures,  OnExtensionModulesOfFiniteHomologicalDimension,FinitenessOfHomologicalDimensionOfExtModules}.

In this paper, we study the consequences of $\RHom$ complexes having finite homological dimension, focusing on projective dimension, Gorenstein projective  dimension and complete intersection dimension. In particular, given $R$-complexes $M$ and $N$, we investigate how the finiteness of projective dimension or complete intersection homological dimension of $\RHom_R(M, N)$ forces the finiteness of the homological dimensions of $M$ or $N$.

We first establish results concerning $\RHom$ complexes with finite complete intersection dimension or finite complete intersection injective dimension; and then derive results concerning Ext modules with finite complete intersection dimension. We address the situation where, given $R$-complexes $M$ and $N$, the derived Hom complex $\RHom_R(M,N)$ has finite Gorenstein projective, projective or complete intersection dimension. We consider two cases: (I) $N$ has finite Gorenstein projective dimension and (II) $M$ has finite Gorenstein injective dimension.

Concerning case (I), under suitable boundedness or finite-generation assumptions, we obtain an isomorphism under non-standard conditions that mimics the tensor evaluation isomorphism:
\begin{theorem}\label{issoss} Let $R$ be a Noetherian ring (not necessarily local) with finite Krull dimension, and let $M$ be a complex in $\mathsf{D}_{\sqsupset}^{\rf}(R)$, $N$ a complex in $\sD(R)$, and $I$ a complex in $\Db(R)$. Assume that $\id_R(I) < \infty$, $\Gpd(\RHom_R(M,N)) < \infty$, and $\Gpd(N) < \infty$. Then the following complexes are isomorphic:
$$\RHom_R(M,N) \ldt_R I \simeq \RHom_R(M,N \ldt_R I).$$
\end{theorem}
As an application of Theorem \ref{issoss}, we investigate the consequences of considering finite projective or complete intersection dimension of the $\RHom$ complex in  case (I), yielding the following stability result:

\begin{theorem}\label{teoinr2}
     Let $M$ and $N$ be non-acyclic complexes in $\Dbf(R)$ such that \linebreak $\RHom_R(M,N)$ is in $\Dbf(R)$.  Assume that  $\gdim_R(N)<\infty$,
     \begin{enumerate}
         \item[(1)] If $\operatorname{pd}_R(\RHom_R(M,N))<\infty$, then $\operatorname{pd}_R(M \oplus N)<\infty$.
         \item[(2)] If $\operatorname{CI-dim}_R(\RHom_R(M,N))<\infty$, then $\operatorname{CI-dim}_R(M \oplus N)<\infty$.
         \end{enumerate}
\end{theorem}

Concerning case (II), under certain boundedness assumption, when we consider finite Gorenstein projective dimension for the corresponding $\RHom$ complex, we obtain an isomorphism that mimic the Hom evaluation isomorphism:
\begin{theorem} \label{HomGidInjs}
Let $R$ be a Noetherian ring (not necessarily local) with finite Krull dimension, and let $M$ and $N$ be complexes in  $\Dlb(R)$.  Assume that $\Gid_R(M)<\infty$ and $\Gpd_R(\RHom_R(M,N))<\infty$. Then for any bounded complex $I$ whose modules are finitely generated and of finite injective dimension, 
  $$I \ldt_R \RHom_R(M,N) \simeq \RHom(\RHom_R(I, M), N).$$
\end{theorem}
As an application of Theorem \ref{HomGidInjs}, we explore the implication of finite complete intersection dimension or finite projective dimension of the $\RHom$ complex in case (II), obtaining:

\begin{theorem}\label{teoinr212}
     Let $M$ and $N$ be non-acyclic complexes in $\Dbf(R)$ such that \linebreak $\RHom_R(M,N)$ is in $\Dbf(R)$.  Assume that $R$ is Cohen-Macaulay and   $\Gid_R(M)<\infty$.
     \begin{enumerate}
         \item  If $\operatorname{pd}_R(\RHom_R(M,N))<\infty$, then $\operatorname{id}_R(M \oplus N)<\infty$.
         \item  If 
$\operatorname{CI-dim}_R(\RHom_R(M,N))<\infty$, then $\operatorname{CI*-id}_R(M \oplus N)<\infty$.
      \end{enumerate}
\end{theorem}

We outline the structure of this paper. In Section 2, we provide some background concerning the derived category of a ring and complete intersection dimensions of complexes. In Section 3, we prove some consequences of having $\RHom$ complexes with finite complete intersection homological dimension and derive some corresponding results for Ext modules. In Section 4, we study the consequences of having $\RHom$ complexes with finite Gorenstein projective, projective or complete intersection dimension, considering the case where the second argument has finite Gorenstein projective dimension. Finally, in Section 5, we consider the case where the first argument has finite Gorenstein injective dimension. 

\section{Background}\label{section:background}

\subsection{The Derived Category of a Ring}
Throughout this paper, we work in the derived category of the ring $R$, denoted $\sD(R)$. This category is obtained from the category of chain complexes by formally inverting all quasi-isomorphisms.
We generally use the same notation and terminology as in this paper as in \cite{DerivedCategoryMethodsInCommutativeAlgebra}. For convenience, we record some definitions and notations in this section.


For a complex $M$ we use $\partial_i^M$ to denote its $i$-th differential $\partial_i^M : M_{i} \to M_{i-1}$. Let $\Dlb (R)$, $\Drb(R)$ and $\Db(R)$ denote the full subcategories of $\sD(R)$ consisting of complexes whose homology is bounded above, bounded below and bounded, respectively. We denote by $\sD^{\rf}(R)$ the full subcategory of $\sD(R)$ consisting of complexes with degreewise finitely generated homology.
The full subcategory $\sD^\rf (R) \cap \Db(R) $ is denoted by $\Dbf (R)$. Similarly, one defines the subcategories $\Dlbf(R)$ and $\Drbf(R)$.

We use the standard derived Hom and tensor product functors. For $R$-complexes $M$ and $N$, we set
\begin{align*}
    \Ext_R^i(M,N)  = \HH_{-i}(\RHom_R(M,N)) \quad \text{and} \quad  \Tor_i^R(M,N)  = \HH_{i} (M \ldt_R N).
\end{align*}

Let $M$ be a complex. The cokernels and kernels of differentials of $M$ are denoted by $C_i(M) = \coker(\partial_{i+1}^M)$ and $Z_i(M) = \ker(\partial_i^M)$. For any integer $n$, the \textit{soft truncations} of $M$ are
\begin{align*}
& M_{\subseteq n}  = 0 \rightarrow C_n(M) \xrightarrow{\overline{\partial}_n^M} M_{n-1} \xrightarrow{\partial_{n-1}^M} M_{n-2} \rightarrow \cdots  \quad  \text{and} \\
&  M_{\supseteq n}  = \cdots \rightarrow M_{n+2} \xrightarrow{\partial_{n+2}^M} M_{n+1} \xrightarrow{\partial_{n+1}^M} Z_n(M) \rightarrow 0,  
\end{align*}
where $\overline{\partial}_n^M: C_n(M) \to M_{n-1}$ is the homomorphism induced by $\partial_n^M$. The definitions of projective dimension Gorenstein projective dimension can be found in \cite{DerivedCategoryMethodsInCommutativeAlgebra}. A fact that we use sometimes without referring is that: If $R$ is a Noetherian ring, not necessarily local, with finite Krull dimension then a complex $M \in \sD(R)$ has finite Gorenstein projective dimension if and only if has finite Gorenstein flat dimension; see \cite[Corollay 17.4.28 and Theorem 19.4.3]{DerivedCategoryMethodsInCommutativeAlgebra}. For more information about the derived category of a ring and its uses in commutative algebra, we recommend Christensen, Foxby, and Holm's excellent and thorough textbook \cite{DerivedCategoryMethodsInCommutativeAlgebra}.

\subsection{Complete Intersection Dimensions} Complete intersection dimensions for complexes were extensively studied by Sather-Wagstaff in \cite{CompleteIntersectionDimensionsForComplexes,CompleteIntersectionDimensionsandFoxbyClasses}.

For a local ring $R$, we remember that a \emph{quasi-deformation} of $R$ is a diagram of local ring homomorphisms  $R \longrightarrow R' \longleftarrow Q$, where $R'$ is a flat $R$-module and where $R'$ is a quotient of $Q$ by a $Q$-regular sequence.


Next, we define complete intersection dimensions in the context of complexes as in \cite{CompleteIntersectionDimensionsandFoxbyClasses} without restricting to homologically bounded complexes. 

\begin{definition}
    Let $M$ be a complex in $\sD(R)$. Then the \emph{complete intersection dimension}, the \emph{complete intersection injective dimension}, and the\emph{upper complete intersection injective dimension} are given by the following formulas:
 {\small 
    \begin{align*}
        \cdim_R(M)&:=\inf \{ \pd_Q(R' \ldt_R M)-\pd_Q(R') \, | \,  R \rightarrow R' \leftarrow Q \mbox{ is a quasi-deformation}\},\\
        \Cid_R(M)&:=\inf \{ \id_Q(R' \ldt_R M)-\pd_Q(R')\, | \,  R \rightarrow R' \leftarrow Q \mbox{ is a quasi-deformation}\}.\\
       \begin{split}
        \uCid_R(M) & := \inf \{ \id_Q(R' \ldt_R M)-\pd_Q(R') \, | \,  R \rightarrow R' \leftarrow Q \text{ is a quasi-deformation} \\ & \hspace{1.3cm}\text{with $R'$ having Gorenstein formal fibers and $R'/\mx R'$ Gorenstein} \}.
\end{split}
         \end{align*}}
       
\end{definition}
When $R$ is a complete intersection, then every bounded complex has finite complete intersection dimension and complete intersection injective dimension. For any  complex $M \in \mathsf{D}(R)$, we immediately have that $\Cid_R(M) \leq \uCid_R(M)$; it is unknown if the two quantities coincide in general. For more information and results concerning complete intersection dimensions, we recommend the following resources \cite{Completeintersectiondimension,CompleteIntersectionDimensionsForComplexes,CompleteIntersectionDimensionsandFoxbyClasses}. 


\begin{notation} Whenever we have a quasi-deformation $R \rightarrow R^\prime \leftarrow Q$, for an $R$-module or complex $M$, we denote $M \otimes_R R^\prime$ by $M^\prime$.
\end{notation}

\begin{remark} \label{rtiq} Given $R$-modules $M$ and $N$, it is unknown if $\cdim_R(M)<\infty$ and $\cdim_R(N)<\infty$ imply that $\cdim_R(M \oplus N)<\infty$. The difficulty of establishing this seemingly simple fact comes from the inability to establish a common quasi-deformation involving $M$ and $N$. As a consequence, we often assume that the direct sum of a collection of modules has finite complete intersection dimension in order to establish a common quasi-deformation. The versions of this open problem with $\Cid$ and $\uCid$ are also unknown. 
 \end{remark}

A result that  will be referenced later in this paper is the following generalization of Ischebeck's formula for finite complete intersection intersection injective dimension. 
\begin{proposition} \label{ischCI-id} Let $M$ and $N$ be non-zero $R$-modules such  that $\Ext_R^i(M,N)=0$ for all $i\gg 0$. Then $$\depth (R)-\depth (M)=\sup\{i\geq 0: \Ext_R^i(M,N)\not=0\},$$
if $\cdim_R (M)<\infty$ or $\Cid_R (N)<\infty$.
\end{proposition}
\begin{proof}
 Assume $\Cid_R (N)<\infty$. By the definition of complete intersection injective dimension, there exists a quasi-deformation $R \rightarrow R' \leftarrow Q$ such that $\id_Q (N')< \infty$. Since $R \rightarrow R'$ is local and flat, we immediately have the following equalities:
    \begin{gather*}\sup\{i \geq 0: \Ext_R^i(M,N) \not=0\}=\sup\{i \geq 0: \Ext_{R'}^i(M',N') \not=0\} \\ \depth_R (R) - \depth_R (M) = \depth_{R'} (R ') - \depth_{R'}( M').\end{gather*} In view of these equalities, we can assume that $R=R'$ and that $R= Q/\boldsymbol{x}Q$ for some $Q$-regular sequence $\boldsymbol{x}=x_1,\dots,x_c$. It then follows that $\sup\{i \geq 0: \Ext_Q^i(M,N) \not=0\}=\depth_Q( Q) - \depth_Q (M)$ by Ischebeck's formula \cite[3.1.24]{bruns}. Moreover, we have that $\sup\{i \geq 0: \Ext_Q^i(M,N) \neq 0 \rbrace =\sup\{i \geq 0: \Ext_R^i(M,N) \neq 0 \rbrace+c$ by \cite[Lemma 2.6]{RemarksonadepthformulaagradeinequalityandaconjectureofAuslander}. Combining the previous equalities, we obtain
\begin{align*}
 \sup\{i \geq 0: \Ext_R^i(M,N) \neq 0 \rbrace  & = (\depth_R (R)+c)-\depth_Q (M) -c
                                              \\ & =\depth_R (R) - \depth_R (M),
 \end{align*}
where the second equality follows using \cite[1.2.26(b)]{bruns}.
\end{proof}
When $\cdim_R (M)<\infty$ instead of $\operatorname{CI-id}_R(N)<\infty$, the same formula holds, this is proved in \cite[Theorem 4.2]{RemarksonadepthformulaagradeinequalityandaconjectureofAuslander}.

To derive some consequences of our results to Ext modules, we frequently use the following lemma.
\begin{lemma}\label{CI-homol}
Let $X$ be a complex in $\Dbf(R)$.  The following results hold:
    \begin{enumerate}
        \item If $\cdim_R(\bigoplus_{i \in \mathbb{Z}} \HH_i(X))<\infty$, then $\cdim_R (X)<\infty$. 
        \item If $\Cid_R(\bigoplus_{i \in \mathbb{Z}} \HH_i(X))<\infty$, then $\Cid_R ( X)<\infty$. 
        \item If $\uCid_R (\bigoplus_{i \in \mathbb{Z}}\HH_i(X))<\infty$, then $\uCid_R (X)<\infty$.
    \end{enumerate}   
\end{lemma}
\begin{proof} We prove item (1). There is a quasi-deformation $R \rightarrow R^\prime \leftarrow Q$ such that $\pd_Q(\bigoplus_{i \in \mathbb{Z}} \HH_i (X^\prime))<\infty$. By \cite[Lemma 2.1(2)]{OnExtensionModulesOfFiniteHomologicalDimension}, we have that $\pd_Q(X^\prime)<\infty$. Then, we have that $\cdim_R (X)<\infty$. 

The proof of items (2) and (3) are similar using \cite[Lemma 4.1]{NumericalApectsofComplexesOfFiniteHomologicalDimensions}.
\end{proof}

\section{Finiteness of CI-dimensions of $\RHom$ and applications to Ext modules}\label{section:FinitenessOfCIdimRHom}

We first work to establish some stability results for $\RHom$ complexes with respect to complete intersection dimension. We need a few necessary lemmas before we start in order to pass homological information between complexes over a ring and over quotients of that ring. The two lemmas below were shared with us by Kaito Kimura.

\begin{lemma}\label{triangle}Let $x \in \mathfrak{m}$ be a non-zero-divisor in $R$, and let $M$ and $N$ be complexes in $\mathsf{D}(R/xR)$.  Then there is a distinguished triangle:   {\small
     $$\RHom_R(M,N) \longrightarrow \RHom_{R/xR}(M,N) \longrightarrow \RHom_{R/xR}(M,N) \longrightarrow \Sigma \RHom_R(M,N)$$ 
     }
\end{lemma}
\begin{proof}
     We first consider the short exact sequence given by multiplication by $x$, namely $$0 \to R \overset{\cdot x}{\longrightarrow} R \longrightarrow R/xR \longrightarrow 0$$ and then apply the functor $\RHom_{R/xR}(M, \RHom_R(-, N))$. By adjunction, this yields the distinguished triangle of the form 
     {\small
     $$\RHom_R(M,N) \longrightarrow \RHom_{R/xR}(M,N) \longrightarrow \RHom_{R/xR}(M,N) \longrightarrow \Sigma \RHom_R(M,N)$$ 
     }
     which is the distinguished triangle we seek.
\end{proof}

\begin{lemma}\label{KimLemma}
 Let $x \in \mx$ be a non-zero-divisor in $R$, and $M$ and $N$ be complexes in $\mathsf{D}(R/xR)$. Suppose that $S$ is a ring such that $R$ is an $S$-algebra. If $\pd_S(\RHom_{R/xR}(M, N))< \infty$, then $\pd_S(\RHom_R(M,N))<\infty$, and if \linebreak $\id_S(\RHom_{R/xR}(M, N))< \infty$, then $\id_S(\RHom_R(M,N))<\infty$.
\end{lemma}
\begin{proof}
This follows from the previous lemma and from the two-out-of-three property for projective and injective dimension \cite[Corollaries 8.1.9 and 8.2.9]{DerivedCategoryMethodsInCommutativeAlgebra}.
\end{proof}

We now obtain the following stability results for complete intersection dimensions.

\begin{proposition}\label{stabilityCI} Let $M$ and $N$ be non-acyclic complexes in $\sD(R)$. For items (1) and (2) below, assume that $M$ and $N$ are in $\Dbf(R)$. For items (3) and (4) assume that $M$ is in $\Drbf (R)$ and $N$
is in $\Dlbf(R)$. For item (4), in addition, assume that $\RHom_R(M,N)$ is in $\Dbf(R)$. The following results hold:
    \begin{enumerate}
        \item If $\cdim_R(\RHom_R(M,N) \oplus M)<\infty$, then $\cdim_R (N)<\infty$.
        \item  If $\cdim_R(\RHom_R(M,N))<\infty$ and $\id_R (N)<\infty$, then $\uCid_R(M)<\infty$.
        \item If $\Cid_R(\RHom_R(M,N))<\infty$, then there is a quasi-deformation for which $\cdim_R (M)$ and $\Cid_R(N)$ are both finite.
        \item If $\uCid_R(\RHom_R(M,N))<\infty$, then there is a quasi-deformation  $R \to R^\prime \leftarrow  Q$ such that $R^\prime/\mathfrak{m}R^\prime $  is artinian Gorenstein, $Q $ is complete and $\pd_Q( M\otimes_R R')<\infty$ and $\id_Q(N\otimes_R R')<\infty$; and therefore $\cdim_R(M) $and $\uCid_R(N)$ are both finite.
    \end{enumerate}
\end{proposition}
\begin{proof}
(1) By the definition of complete intersection dimension, there is a quasi-deformation $R \rightarrow R^\prime \leftarrow Q$ such that $\pd_Q(\RHom_{R^\prime}(M^\prime, N^\prime))<\infty$ and $\pd_Q (M^\prime)<\infty$. By Lemma \ref{KimLemma}, $\pd_Q(\RHom_{Q}(M^\prime, N^\prime))<\infty$. Now, by \cite[Proposition 16.4.19]{DerivedCategoryMethodsInCommutativeAlgebra}, it follows that $\pd_Q (N^\prime)<\infty$. Hence, $\cdim_R(N)<\infty$.

(2) Note that $\RHom_R(M,N)$ is homologically bounded since $\operatorname{id}_R(N)<\infty$. By \cite[Theorem F]{CompleteIntersectionDimensionsandFoxbyClasses}, there is a quasi-deformation $R \rightarrow R^\prime \leftarrow Q$ with $R^\prime/\mathfrak{m}R^\prime$ artinian Gorenstein and $Q$ complete having the property that $\pd_Q(\RHom_{R^\prime}(M^\prime, N^\prime))<\infty$. By Lemma \ref{KimLemma}, $\pd_Q(\RHom_{Q}(M^\prime, N^\prime))<\infty$. On the other hand, since $\id_R (N)<\infty$, we have that by \cite[Theorem 16.4.36]{DerivedCategoryMethodsInCommutativeAlgebra}, $\id_{R^\prime} (N^\prime)<\infty$. From \cite[Proposition 15.4.10]{DerivedCategoryMethodsInCommutativeAlgebra} we have that $\id_Q (N^\prime)<\infty$ as $R'$ has finite projective dimension over $Q$. Hence $\id_Q(M^\prime)<\infty$ by \cite[Proposition 16.4.34]{DerivedCategoryMethodsInCommutativeAlgebra}, and so it follows that $\uCid_R (M)<\infty$ by the definition of upper CI-injective dimension.

(3) There is a quasi-deformation $R \rightarrow R^\prime \leftarrow Q$ such that $\id_Q(\RHom_{R^\prime}(M^\prime, N^\prime))<\infty$. By Lemma \ref{KimLemma}, $\id_Q(\RHom_Q(M^\prime, N^\prime))<\infty$. Then by \cite[Proposition 16.4.32]{DerivedCategoryMethodsInCommutativeAlgebra}, $\pd_Q (M^\prime)<\infty$ and $\id_Q (N^\prime)<\infty$. Thus, there is a quasi-deformation that gives the finiteness of $\cdim_R (M)$ and $\Cid_R(N)$.

(4) This is similar to the proof of the previous item using \cite[Proposition 3.5]{CompleteIntersectionDimensionsandFoxbyClasses}.
\end{proof}
Motivated by the results in \cite{GorensteinringsviahomologicaldimensionsandsymmetryinvanishingofExtandTatecohomology,NumericalApectsofComplexesOfFiniteHomologicalDimensions,OnExtensionModulesOfFiniteHomologicalDimension,FinitenessOfHomologicalDimensionOfExtModules}, which investigate what finite homological dimensions of Ext modules reveal about the original modules, we apply our results to address analogous questions for complete intersection dimensions.

\begin{corollary}
    \label{teo1.2}  Let $M$ and $N$ be non-zero $R$-modules such that $\Ext_R^i(M,N)=0$ for all $i\gg 0$. Set $L=\bigoplus_{i=0}^\infty \Ext_R^i(M,N)$.  The following results hold:
 \begin{enumerate}
\item If $\cdim_R(M \oplus L )<\infty$, then  $\cdim_R(N)<\infty.$
\item  If $\cdim_R(L)<\infty$ and $\id_R (N)<\infty$, then $\uCid_R(M)<\infty$.
 \item  If $\Cid_R(L)<\infty$, then there is a quasi-deformation for which $\cdim_R(M)$ and $\Cid_R(N)$ are both finite.
 \item If $\uCid_R(L)<\infty$, then there is a quasi-deformation for which $\cdim_R(M)$ and $\uCid_R(N)$ are both finite.
 \end{enumerate}
 \end{corollary}
 \begin{proof}
     (1)  Applying Lemma \ref{CI-homol}(1) to the complex $\RHom_R(M,N) \oplus M$, we have that $\cdim_R (\RHom_R(M,N) \oplus M)<\infty$. Thus, by Proposition \ref{stabilityCI}(1), \linebreak $\cdim_R (N)<\infty$. 

     (2) By Lemma \ref{CI-homol}(1), we have that $\operatorname{CI-dim}_R(\RHom_R(M,N))<\infty$. Hence, the conclusion follows from Proposition ~\ref{stabilityCI}(2).



     (3) Lemma \ref{CI-homol}(2) shows that $\operatorname{CI-id}_R(\RHom_R(M,N))<\infty$. Thus follows from Proposition \ref{stabilityCI}(3).

     (4) This is similar to the previous item but uses Lemma \ref{CI-homol}(3) and Proposition \ref{stabilityCI}(4).
 \end{proof}

With the notation of Corollary \ref{teo1.2}, we give below an example where
$\operatorname{CI-dim}_R(L)<\infty$, over a ring that is not necessarily a complete intersection. 

\begin{example}
Let $S[[x,y]]$ be a formal power series ring in two variables over the local ring $S$. Consider the following ring and the following modules: $R=S[[x,y]]/(x^2,y^2)$, $M=R/(x)$ and $N=R/(y)$. We then claim that:
\begin{enumerate}
    \item $\operatorname{CI-dim}_R(M)=\operatorname{CI-dim}_R(N)=0$.
    \item $\Ext_R^i(M,N)=0$ for all $i>0$ and $\pd_R(M)=\pd_R(N)=\infty$.
    \item $\cdim_R(\Hom_R(M,N))=0$.
\end{enumerate}
\begin{proof}
(1) Setting $Q= S[[x,y]]/(y^2)$, we may consider a trivial quasi-deformation $R \rightarrow R \leftarrow Q$. Note that $M=R/(x) \cong Q/(x)$. Since $\pd_Q(M)- \pd_Q(R)=1-1$, we have $\cdim_R(M)=0$. A similar argument shows that $\cdim_R(N)=0$.

(2) Consider the minimal free resolution: 
\begin{align*}
    F= \cdots \rightarrow R \xrightarrow{x} R \xrightarrow{x} R \xrightarrow{x} R \rightarrow 0
\end{align*}
of $M$ over $R$. We have that $\pd_R(M)= \infty$. The complex $\Hom_R(F,N)$ is 
\begin{align*}
    0 \rightarrow N \xrightarrow{x} N \xrightarrow{x} N\xrightarrow{x} N \rightarrow \cdots.
\end{align*}
Thus, computing homologies we see that $\Ext_R^i(M,N)=xN/xN=0$ for all $i>0$. We remark that a similar argument shows that $\pd_R(N)= \infty$.

(3) Note that $\Hom_R(M,N) \cong xN \cong N/xN \cong R/(x,y)$. Setting $T=S[[x,y]]$, we may consider a trivial quasi-deformation $R \rightarrow R \leftarrow T$. Observe that $\Hom_R(M,N) \cong T/(x,y)$ and that $\operatorname{pd}_T(R)=2$. Thus, $\pd_T(\Hom_R(M,N))-\pd_T(R)=2-2=0$ and $\operatorname{CI-dim}_R(\Hom(M,N))=0$.
\end{proof}
\end{example}

To complement \cite[Theorem 1.1]{OnExtensionModulesOfFiniteHomologicalDimension} by including complete intersection dimension and upper complete intersection injective dimension, we derive the following corollary: 

\begin{corollary}
    Let $M$ and $N$ be non-zero $R$-modules, and let $r=\operatorname{Rfd}_R (M)$.  The following results hold: 
    \begin{enumerate}
        \item Assume that $\operatorname{CI-dim}_R ( \bigoplus_{i=0}^r\Ext_R^i(M,N))<\infty$.
        \begin{enumerate}
            \item[(1.a)] If  $\operatorname{pd}_R (M)<\infty$, then $\operatorname{CI-dim}_R(N)<\infty$.
            \item[(1.b)] If  $\operatorname{id}_R (N)<\infty$,  then $\operatorname{CI*-id}_R(M)<\infty$.  
        \end{enumerate}
        \item Assume that $\operatorname{CI*-id}_R (\bigoplus_{i=0}^r \Ext_R^i(M,N))<\infty$. 
       \begin{enumerate}
           \item[(2.a)]  If $\operatorname{pd}_R (M)<\infty$, then $\operatorname{CI*-id}_R(N)<\infty$. 
            \item[(2.b)]   If $\id_R(N)<\infty$, then $\operatorname{CI-dim}_R(M)<\infty$.
       \end{enumerate}
       \item Assume that $\operatorname{CI-id}_R (\bigoplus_{i=0}^r \Ext_R^i(M,N))<\infty$.
       \begin{enumerate}
           \item[(3.a)]  If $\operatorname{pd}_R (M)<\infty$, then $\operatorname{CI-id}_R(N)<\infty$. 
            \item[(3.b)]   If $\id_R(N)<\infty$, then $\operatorname{CI-dim}_R(M)<\infty$.
       \end{enumerate}
       
\end{enumerate}
\end{corollary}
\begin{proof} If $\operatorname{pd}_R (M)<\infty$, we have that $r=\operatorname{pd}_R (M)$ and $\Ext_R^i(M,N)=0$ for all $i>r$.  If $R$ is Cohen-Macaulay, then $r=\depth(R)-\depth(M)$,  and  moreover if $N$ has finite injective dimension, then $\Ext_R^i(M,N)=0$ for all $i>r$ by Ischebeck's formula. Thus, note that in each one of the sub-items, one has that $\Ext_R^i(M,N)=0$ for all $i>r$. Hence, the corollary follows from the previous theorem. 
\end{proof}

Let $M$ and $N$ be non-zero $R$-modules. The equality
$$\operatorname{Rfd}_R M=\depth(R)-\depth(M)=\sup\{i \geq 0\mid ~ \Ext_R^i(M,N)\not=0\}$$
holds in certain cases provided that $\Ext_R^i(M,N)=0$ for all $i\gg 0$. For instance, this holds if  $\operatorname{CI-dim}_R(M)<\infty$ (\cite[Theorem 4.2]{RemarksonadepthformulaagradeinequalityandaconjectureofAuslander}) or $\operatorname{CI-id}_R(N)<\infty$ (Proposition \ref{ischCI-id}). One can use this together with Corollary \ref{teo1.2} to obtain results on the finiteness of the complete intersection dimension of $\bigoplus_{i=0}^{\operatorname{Rfd}_R (M)} \Ext_R^i(M,N)$ provided that $\Ext_R^i(M,N)=0$ for all $i\gg 0$. Namely, we have the following corollary:

 \begin{corollary}\label{cor3.7} Let $M$ and $N$ be non-zero $R$-modules such that $\Ext_R^i(M,N)=0$ for all $i\gg 0$. Set $L=\bigoplus_{i=0}^{\Rfd_R(M)} \Ext_R^i(M,N)$.  The following results hold:
\begin{enumerate}
    \item If $\operatorname{CI-dim}_R(M \oplus L)<\infty$, then $\operatorname{CI-dim}_R(N)<\infty$.
    \item If $\operatorname{CI-id}_R(L)<\infty$  and $\operatorname{CI-id}_R(N)<\infty$, then $\operatorname{CI-dim}_R(M)<\infty$. 
\end{enumerate}
\end{corollary}

We finish this section by studying deficiency modules with finite complete intersection homological dimension, all of them with the same quasi-deformation. This is motivated by previous work where the consequences of deficiency modules having other finite homological dimensions were studied: In \cite{FH}, the authors investigated the consequences of the finiteness of projective and injective dimensions of the deficiency modules; while in \cite[Theorem 3.1]{OnExtensionModulesOfFiniteHomologicalDimension} this direction was considered with Gorenstein dimensions. 

For the convenience of the reader, we recall the definition of deficiency modules \cite{Sch98}. If $R$ is a local ring with a normalized dualizing complex $D$, then for an $R$-module $M$ and for each $i$, the $R$-module  
$$ K^i(M): = \operatorname{H}_i(\RHom_R(M,D))$$
is called the $i$-th {\it deficiency module} of $M$. It is known that $K^i(M)=0$ for $i< \depth(M)$ or $i> \dim (M)$, and $K^i(M)\neq 0$ for $i = \depth(M)$ and $i= \dim(M)$. For more about deficiency modules, we recommend the reference \cite{Sch98}.



Now,  we  provide a version of \cite[Theorem 3.1]{OnExtensionModulesOfFiniteHomologicalDimension} with complete intersection dimension and  complete intersection injective dimensions. 
\begin{theorem}\label{defic1}
    Let $R$ have a normalized dualizing complex $D$, and let $M$ be a non-zero $R$-module.  The following results hold:
    \begin{enumerate}
        \item If $\cdim_R\left(\bigoplus\limits_{i = \depth(M)}^{\dim(M)}K^i(M)\right)< \infty$, then $\uCid_R(M)<\infty$ and $R$ is Cohen-Macaulay.
        \item If $\Cid_R\left(\bigoplus\limits_{i = \depth(M)}^{\dim(M)}K^i(M)\right)< \infty$, then $\cdim_R(M)<\infty$ and $R$ is Cohen-Macaulay.
    \end{enumerate}  
\end{theorem}
\begin{proof}
(1) First, observe that $\RHom_R(M,D) \in \Dbf(R)$. Moreover, using \linebreak Lemma~\ref{CI-homol}(1) one can see that $\cdim_R(\RHom_R(M, D)) < \infty$, since the modules of deficiency are simply the homology modules of $\RHom_R(M, D)$. We then have that $\uCid_R(M) < \infty$ by Proposition \ref{stabilityCI}(2) or  \cite[Corollary 4.6(b)]{CompleteIntersectionDimensionsandFoxbyClasses}. This implies that $R$ is Cohen-Macaulay by \cite[Theorem 2.7]{CharacterizingLocalRingsViaCompleteIntersectionHomologicalDimensions}.

(2) Using Lemma~\ref{CI-homol}(2), one can see that $\Cid_R(\RHom_R(M,D))< \infty$. By Proposition \ref{stabilityCI}(3), we have $\cdim_R(M)< \infty$, as required.

\end{proof}

\section{Finite homological dimensions of $\RHom$ with finite Gpd in the second argument}

In this section, we study the consequences of having $\RHom$ complexes of finite Gorenstein projective, projective or complete intersection dimension where its second argument has finite Gorenstein projective dimension. 

To prepare the mains result of this section, we  first prove a necessary acyclicity lemma first.  This lemma will also be need to prove one of the main theorems in the next section.

\begin{lemma}\label{seqext} Let $R$ be a Noetherian ring with finite Krull dimension. Let $(X, \partial)$ be an acyclic complex of Gorenstein flat $R$-modules. If $\coker (\partial_n)$ has finite Gorenstein flat dimension for some $n$, then  $\coker (\partial_j)$ is  a Gorenstein flat module for each $j$; and if in addition $I$ is a module of finite injective dimension, then $X \otimes_R I$ is acyclic.
\end{lemma}
\begin{proof}Since $X$ is acyclic, for any $i, j \in \mathbb{Z}$ such that $j \geq i$ we have an exact diagram    
{\tiny
\begin{center}
\begin{tikzcd}[row sep = small, column sep = small] 
\cdots \ar[r] & 
X_j \ar[r, "\partial_j"] & 
X_{j-1} \ar[rr, "\partial_{j-1}"] \ar[dr] && 
X_{j-2} \ar[r] & 
\cdots \ar[r] & 
X_{i-1} \ar[rr, "\partial_i"] \ar[rd] &  &X_{i-2} \ar[r]& \cdots  \\
&&& \coker(\partial_{j}) \ar[ur] \ar[rd] &&&& \coker(\partial_{i}) \ar[ru] \ar[dr] \\
&& 0 \ar[ru] && 0&&0 \ar[ru] && 0
\end{tikzcd}
\end{center}}
From the the two-out-of-three property for Gorenstein flat dimension (\cite[Proposition 9.3.25]{DerivedCategoryMethodsInCommutativeAlgebra}), note that if one of $\coker(\partial_i)$ or $\coker(\partial_j)$ has finite Gorenstein flat dimension, then the other also has finite Gorenstein flat dimension. Thus, $\coker(\partial_j)$ has finite Gorenstein dimension for all $j$. 

Set $t=\sup \{ \depth (R_\mathfrak{p}): \mathfrak{p} \in \Spec(R)\}$. Since $R$ has finite Krull dimension, we have that $t$ is finite. Now, fix $j \in \mathbb{Z}$ and take $i$ such that $j>t+i+1$, in the diagram. We then obtain an exact sequence
$$0 \rightarrow \coker(\partial_j) \rightarrow  X_{j-2} \rightarrow \cdots \rightarrow X_{i}.$$
By repeated applications of the depth lemma \cite[Theorem 14.3.20]{DerivedCategoryMethodsInCommutativeAlgebra}, we have that $\depth_{R_\px}(\coker(\partial_j)_\px) \geq \depth_{R_\px}(R_\px)$ for all $\px \in \Spec(R)$. It then follows from \cite[Theorem 19.3.6]{DerivedCategoryMethodsInCommutativeAlgebra} that $\coker( \partial_j)$ is Gorenstein flat. Thus, the first assertion is proved.

We now prove that $X \otimes I$ is acyclic. We immediately have that $\Tor^R_i(I, X_j) = 0$ for all $i > 0$ and $j \in \mathbb{Z}$ by \cite[Theorem 9.3.41]{DerivedCategoryMethodsInCommutativeAlgebra}. Since $\coker(\partial_j)$ is Gorenstein flat, we have that $\Tor^R_i(I, \coker(\partial_j)) = 0$ for all $i > 0$ again by \cite[Theorem 9.3.41]{DerivedCategoryMethodsInCommutativeAlgebra}. Thus, by \cite[Proposition A.10]{DerivedCategoryMethodsInCommutativeAlgebra}, we have that $X \otimes_R I$ is acyclic.
\end{proof}
Now, we prove our first main result, which 
establishes a type of tensor-evaluation isomorphism in the derived category under non-standard assumptions.
\begin{theorem}\label{issos} Let $R$ be a Noetherian ring (not necessarily local) with finite Krull dimension, and let $M$ be a complex in $\mathsf{D}_{\sqsupset}^{\rf}(R)$, $N$ a complex in $\sD(R)$, and $I$ a complex in $\Db(R)$. Assume that $\id_R(I) < \infty$, $\Gpd(\RHom_R(M,N)) < \infty$, and $\Gpd(N) < \infty$. Then the following complexes are isomorphic:
$$\RHom_R(M,N) \ldt_R I \simeq \RHom_R(M,N \ldt_R I).$$
\end{theorem}
\begin{proof} 
Let $F$ and $G$ be semi-projective replacements of $M$ and $N$, respectively, and let $J$ be a semi-injective replacement of $I$. We may assume that $F$ is a complex of finitely generated projective modules bounded below and that $J$ is bounded. Set $n=\Gfd_R (N)$. By \cite[9.3.19]{DerivedCategoryMethodsInCommutativeAlgebra}, $n\geq \sup G$ and therefore the natural morphism $\tau_{\subseteq n}: G \to G_{\subseteq n}$ is a quasi-isomorphism. Since $F$ is semi-projective, then $\Hom_R(F, \tau_{\subseteq n}): \Hom_R(F, G) \to \Hom_R(F,G_{\subseteq n})$ is a quasi-isomorphism. Let $\Theta: P \to \Hom_R(F, G)$ be a semi-projective resolution of $\Hom_R(F, G)$, and consider the morphism $\Phi = \Hom_R(F, \tau_{\subseteq n}) \circ \Theta$, which is a quasi-isomorphism. We will show that $\Phi \otimes J$ is a quasi-isomorphism.

Let $L = \Hom_R(F, G_{\subseteq n})$, and note that $L$ is a bounded above complex of Gorenstein flat modules. Let $C=\cone(\Phi)$. We have that $C_i= L_i \oplus P_{i-1}$ for each $i$, and so $C$ is an acyclic complex of Gorenstein flat modules. Let $m=\Gfd_R (\RHom_R(M,N))$. For $i > \sup(L^\natural)$, we have that $\CC_{i-1}(P)=\CC_i(C)$ (where $\CC_i$ represents the cokernel at index $i$.) Since $\CC_i(P)$ is Gorenstein flat for $i > \Gfd_R(\RHom_R(M, N))$ by \cite[Lemma 9.3.26]{DerivedCategoryMethodsInCommutativeAlgebra}, we have that the cokernels of $C$ are Gorenstein flat for sufficiently large indexes. We then have that $C \otimes_R J_v$ is acyclic for all $v \in \mathbb{Z}$ by Lemma~\ref{seqext}. Since $J$ is bounded below, we have that $C \otimes_R J$ is acyclic by \cite[Proposition A.11]{DerivedCategoryMethodsInCommutativeAlgebra}. This implies that $\Phi \otimes J$ is a quasi-isomorphism.

We then have that $P \otimes_R J$ is quasi-isomorphic to $\Hom_R(F, G_{\subseteq n}) \otimes_R J$. But by \cite[Theorem 4.5.10(1)-(a)]{DerivedCategoryMethodsInCommutativeAlgebra}, we have the quasi-isomorphism $\Hom_R(F, G_{\subseteq n}) \otimes_R J \simeq \Hom_R(F, G_{\subseteq n} \otimes_R J)$. Thus,
\begin{equation*}
   P \otimes_R J \simeq \Hom_R(F, G_{\subseteq n} \otimes_R J).
\end{equation*}
But $P \otimes_R J \simeq \RHom_R(M, N) \ldt_R I$ and $\Hom_R(F, G_{\subseteq n} \otimes_R J) \simeq \RHom_R(M, N \ldt_R I)$ in $\sD(R)$ (noting that the last isomorphism uses \cite[Corollary 9.3.30]{DerivedCategoryMethodsInCommutativeAlgebra}). We then have the desired isomorphism
\begin{equation*}
\RHom_R(M,N) \ldt_R I \simeq \RHom_R(M,N \ldt_R I).
\end{equation*}
\end{proof}

Now we prove our second main theorem, where $\hdim$ denotes projective dimension or complete intersection dimension.

\begin{theorem}\label{reqew} 
    Let $M$ and $N$ be non-acyclic complexes in $\Dbf(R)$  such that \linebreak $\RHom_R(M,N)$ is in $\Dbf(R)$. Assume that $\gdim_R(N)<\infty$. 
       If \linebreak $\hdim_R(\RHom_R(M,N))<\infty$ , then  $\hdim_R(M \oplus N)<\infty$. 
\end{theorem}
\begin{proof}
We may assume that $R$ is complete, and so has a dualizing complex $D$.  
In particular, \linebreak $\gdim_R(\RHom_R(M,N))<\infty$. Then, by Theorem \ref{issos}, 
\begin{equation}\label{opoes}
    \RHom_R(M, N \ldt_R D)\simeq  \RHom_R(M,N) \ldt_R D.
\end{equation}

If $\hdim_R=\pd_R$, then  by  \eqref{opoes} and Foxby-Sharp equivalence \cite[Theorem 10.3.8]{DerivedCategoryMethodsInCommutativeAlgebra}, $\RHom_R(M, N \ldt_R D)$ is homologically bounded and of finite injective dimension. Now, by \cite[Proposition 16.4.32]{DerivedCategoryMethodsInCommutativeAlgebra}, $\pd_R (M)<\infty$ and $\id_R (N \ldt_R D)<\infty$. The last inequality implies that $\pd_R (N)<\infty$ through Foxby-Sharp equivalence \cite[Theorem 10.3.8]{DerivedCategoryMethodsInCommutativeAlgebra}.

If $\hdim_R=\cdim_R$, then by \eqref{opoes} and \cite[Theorem E(a)]{CompleteIntersectionDimensionsandFoxbyClasses}, we have that $\RHom_R(M, N \ldt_R D)$ is homologically bounded and  $\uCid_R(\RHom_R(M, N \ldt_R D))<\infty$. We then have, by Proposition \ref{stabilityCI}(4), a quasi-deformation $R \to R^\prime \leftarrow Q$ such that $R^\prime/\mathfrak{m}R^\prime$ is artinian Gorenstein, $\pd_Q(M \otimes_R R^\prime)<\infty$, and $\id_Q( (N \ldt_R D) \ldt_R R^\prime)<\infty$. Since $R^\prime/\mathfrak{m}R^\prime$ is Gorenstein, then  $R^{\prime} \ldt_R D$ is dualizing for $R^{\prime}$. As $Q$ is complete, it admits a dualizing complex $D^Q$. As $Q\to R'$ is surjective with kernel generated by a $Q$-regular sequence, we have that $R^{\prime} \ldt_Q D^Q$ is dualizing for $R^{\prime}$. We then have that $R^{\prime} \ldt_Q D^Q \simeq R^{\prime} \ldt_R D$ (up to a shift). This yields the following isomorphisms (up to a shift)
\begin{align*}
    (N \ldt_R D) \ldt_R R^\prime    &\simeq N \ldt_R (D \ldt_R R^\prime)\\
                                    &\simeq N \ldt_R (R^\prime \ldt_Q  D^Q)\\
                                    &\simeq (N \otimes_R R^\prime) \ldt_Q D^Q.
\end{align*}
Hence, $\id_Q( (N \otimes_R R^\prime) \ldt_Q D^Q)<\infty$, so $\pd_Q(N \otimes_R R^\prime)<\infty$ through Foxby-Sharp equivalence \cite[Theorem 10.3.8]{DerivedCategoryMethodsInCommutativeAlgebra}. Consequently, $\cdim_R(M \oplus N)<\infty$. 
\end{proof}

\begin{corollary} Let $M$ and $N$ be non-acyclic complexes in $\Dbf(R)$. Assume that $\pd_R(\RHom_R(M,N))<\infty$ and that $\RHom_R(M,N)$ is in $\Dbf(R)$. Then the following are simultaneously finite: $\pd_R(M)$, $\cdim_R(M)$, $\pd_R(N)$, and $\gdim_R(N)$. 
\end{corollary}
\begin{proof} If $\pd_R(M)<\infty$ (resp. $\pd_R(N)<\infty$), then $\cdim_R(M)<\infty$ (resp. $\gdim_R(N)<\infty$). Thus, if we show the implications below, we will be done. 
$$\xymatrix{\cdim _R (M) < \infty \ar@{=>}[r]^{(1)} & \gdim_R (N)< \infty \ar@{=>}[r]^{(2)} & \pd_R (M \oplus N) < \infty  }$$

To prove (1), assume that $\cdim_R(M)<\infty$. Then $\cdim_R(\RHom_R(M,N) \oplus M)<\infty$, by \cite[Lemma 3.6]{CompleteIntersectionDimensionsForComplexes}. Thus,  by Proposition \ref{stabilityCI}(1), $\cdim_R(N)<\infty$ and then $\gdim_R (N) < \infty$.

To prove (2), note that if $\gdim_R(N)<\infty$, then by Theorem \ref{reqew}, $\pd_R(M \oplus N)<\infty$.
\end{proof}

The following example shows that the condition $\pd_R(\RHom_R(M,N)) < \infty$ in the above corollary cannot be replaced by the weaker assumption \linebreak $
\cdim_R(\RHom_R(M,N)) < \infty$.

\begin{example}
Let $S = T[[x]]$ be the formal power series ring in one variable over a local ring $T$. Set $R := S/(x^n)$, where $n \ge 2$, and let $M := Rx$. Then $\cdim_R(\RHom_R(M,R)) = 0$ and $\cdim_R(M) = 0$, while $\pd_R(M) = \infty$; see \cite[Example~3.11]{ARCForModulesWhose(Self)DualHasFiniteCompleteIntersectionDimensionv2}.
\end{example}

\begin{corollary} \label{fint} Let $M$ be a non-zero complex in $\Dbf(R)$ such that $\RHom_R(M,M)$ is in $\Dbf(R)$. Assume that $\cdim_R(\RHom_R(M,M))<\infty$. If  $\gdim_R(M)<\infty$, then $\pd_R (M)<\infty$.
\end{corollary}
\begin{proof} By Theorem \ref{reqew}, $\cdim_R(M)<\infty$. Then, by \cite[Proposition 5.7(1)]{PersistenceOfHomologyOverCommutativeNoetherianRings}, $\pd_R(M)<\infty$.
\end{proof}
\section{Finite homological dimension of $\RHom$ with finite Gid in the first argument}

In this section, we study the consequences of having $\RHom$ complexes of finite Gorenstein projective, projective or complete intersection dimension where its first argument has finite Gorenstein injective dimension. Next, our main first result in this section.

\begin{theorem} \label{HomGidInj}
Let $R$ be a Noetherian ring (not necessarily local) with finite Krull dimension, and let $M$ and $N$ be complexes in  $\Dlb(R)$.  Assume that $\Gid_R(M)<\infty$ and $\Gpd_R(\RHom_R(M,N))<\infty$. Then for any bounded complex $I$ whose modules are finitely generated and of finite injective dimension, 
  $$I \ldt_R \RHom_R(M,N) \simeq \RHom(\RHom_R(I, M), N).$$
\end{theorem}
\begin{proof}
Let $E$ and $J$ be semi-injective replacements of $M$ and $N$ respectively. By assumption, we may assume that $E$ and $J$ are both bounded above. Set $n=\Gid_R (M)$. Then the natural map $\tau_{\supseteq -n}: E_{\supseteq -n} \to E$ is a quasi-isomorphism as $n \geq -\inf(M)$. Since $J$ is semi-injective, we have that $\Hom_R(\tau_{\supseteq -n}, J)$ is a quasi-isomorphism. Let $\Theta: P \to \Hom_R(E, J)$ be a semi-projective resolution, and consider the composition $\Phi: \Hom_R(\tau_{\supseteq -n}, J) \circ \Theta$. We will prove that $I \otimes \Phi$ is a quasi-isomorphism.
   
Set $L=\Hom_R(E_{\supseteq - n}, J)$ and $C=\cone(\Phi)$. Note that $L$ is a bounded above complex of Gorenstein flat modules. Since $C_i=  L_{i} \oplus P_{i-1}$ for each $i$, we have that $C$ is an acyclic complex of Gorenstein flat modules. Let $m=\Gfd_R (\RHom_R(M,N))$. By \cite[Lemma 9.3.26]{DerivedCategoryMethodsInCommutativeAlgebra}, we have that $\CC_i(P)$ is Gorenstein flat for all $i \geq m$, which forces $\CC_j(C)$ to be Gorenstein flat for $j \gg 0$ as $L$ is bounded above. Since each module $I_v$ of $I$ has finite injective dimension, Lemma \ref{seqext} implies that $I_v \otimes_R C$ is acyclic for each $v \in \mathbb{Z}$, and so by \cite[Proposition A.11]{DerivedCategoryMethodsInCommutativeAlgebra}, $I \otimes_R C$ is acyclic. 

We now have that $I \otimes_R P$ is quasi-isomorphic to $I \otimes_R \Hom_R(E_{\supseteq -n}, J)$. By Hom-evaluation (\cite[Theorem 12.1.16(d)]{DerivedCategoryMethodsInCommutativeAlgebra}) we have that $I \otimes_R \Hom_R(E_{\supseteq -n}, J)$ is quasi-isomorphic to $\Hom_R(\Hom_R(I, E_{\supseteq -n}), J)$. 
In addition, as $E_{\supseteq -n}$ is a bounded complex of Gorenstein injective modules, by  \cite[Lemma 3.4]{TheAuslanderBoundForComplexes} we obtain that \linebreak $\Hom_R(I, E_{\supseteq -n}) \simeq \Hom_R(I, E)$  once that is noticed that the Ext vanishings needed in that lemma are a result of \cite[Theorem 9.2.24]{DerivedCategoryMethodsInCommutativeAlgebra}. 

Thus, we have that $I \otimes_R P$ is quasi-isomorphic to $\Hom_R(\Hom_R(I, E), J)$, which gives the desired isomorphism in $\sD(R)$:
\begin{equation*}
   I \ldt_R \RHom_R(M,N) \simeq \RHom_R(\RHom_R(I, M), N).
\end{equation*}

\end{proof}

As an application of Theorem \ref{HomGidInj}, we derive the next stability result which is the second main result of this section. It is dual to Theorem~\ref{reqew}, and the notation used is the following: $\hdim$ denotes projective dimension or complete intersection case, and $\hid$ denotes  injective dimension in the first case, while in the second case denotes upper complete intersection injective dimension.

\begin{theorem}\label{reqs22} Let $R$ be a Cohen-Macaulay local ring, and let $M$ and $N$ be non-acyclic complexes in $\Dbf(R)$ such that $\RHom_R(M,N)$ is in $\Dbf(R)$. Assume that $\Gid_R(M)<\infty$. If $\hdim_R(\RHom_R(M,N))<\infty$, then $\hid_R(M \oplus N)<\infty$. 
\end{theorem}
\begin{proof}
We may assume that $R$ is complete, and therefore $R$ admits a canonical module $\omega$. By assumption, we have $\gdim_R(\RHom_R(M,N))<\infty$. Therefore by Theorem \ref{HomGidInj}, 
\begin{equation}\label{fetw}
    \RHom_R(\RHom_R(\omega, M), N) \simeq \omega \ldt_R \RHom_R(M,N).
\end{equation}

If $\hdim_R=\pd_R$, then we have that $\RHom_R(\RHom_R(\omega, M),N)$ is homologically bounded of finite injective dimension by \eqref{fetw} and Foxby-Sharp equivalence. By Bass series arguments \cite[Proposition 16.4.32]{DerivedCategoryMethodsInCommutativeAlgebra}, this implies that \linebreak $\pd_R(\RHom_R(\omega, M)) < \infty$ and $\id_R(N) < \infty$. Again through Foxby-Sharp equivalence, we have that $\id_R(M)<\infty$.

Suppose now that $\hdim_R=\cdim_R$. Then we have that  \linebreak $\RHom_R(\RHom_R(\omega, M), N)$ is homologically bounded and \linebreak $\uCid_R(\RHom_R(\RHom_R(\omega, M), N)) < \infty$ by \eqref{fetw} and Foxby-Sharp equivalence. By Proposition \ref{stabilityCI}(4), there is a quasi-deformation $R \to R^\prime \leftarrow Q$ such that $R^\prime/\mathfrak{m}R^\prime$ is artinian Gorenstein, $Q$ is complete, $\pd_Q(\RHom_R(\omega,M) \otimes_R R^\prime )<\infty$, and $\id_Q(N \otimes_R R^\prime)<\infty$. Since $R'$ is flat over $R$, we have that $\RHom_R(\omega, M) \otimes_R R' \simeq \RHom_{R'}(\omega \otimes_R R', M \otimes_R R')$. Note that $Q$ is Cohen-Macaulay with a canonical module $\omega_Q$ since it is complete. In addition, note that $\omega^\prime = \omega \otimes_R R^\prime$ is a canonical module for $R^\prime$ since $R'/\mx R'$ is Gorenstein, and that $ \omega^\prime \cong \omega_Q \otimes_Q R^\prime \simeq \omega_Q \ldt_Q R^\prime$. Since $$\RHom_{R'}(\omega', M \otimes_R R') \simeq \RHom_{R'}(\omega_Q \ldt_Q R', M \otimes_R R') \simeq \RHom_Q(\omega_Q, M \otimes_R R'),$$ we have that $\pd_Q(\RHom_Q(\omega_Q, M \otimes_R R'))<\infty$, whence it follows that $\id_Q(M \otimes_R R')<\infty$ using Foxby-Sharp equivalence. Thus, $\uCid_R(M \oplus N)<\infty$.
\end{proof}

\begin{corollary} \label{fint2} Let $R$ be a Cohen-Macaulay local ring, and let $M$ be a non-acyclic complex in $\Dbf(R)$ such that $\RHom_R(M,M)$ is in $\Dbf(R)$. Assume that $\cdim_R(\RHom_R(M,M))<\infty$.  If $\Gid_R(M)<\infty$, then $\id_R(M)<\infty$.
\end{corollary}
\begin{proof} By Theorem \ref{reqs22}, $\uCid_R(M)<\infty$. Then, by \cite[Proposition 5.7(2)]{PersistenceOfHomologyOverCommutativeNoetherianRings}, $\id_R(M)<\infty$.
\end{proof}

The following corollary is a version of \cite[Theorem 1.2]{FinitenessOfHomologicalDimensionOfExtModules} for complete intersection dimension.

\begin{corollary}
    Let $M$ and $N$ be non-zero $R$-modules, and assume that \linebreak $\cdim_R(\bigoplus_{i=0}^{\Rfd_R (M)} \Ext_R^i(M,N))<\infty$. If $\id_R(M)<\infty$, then $\uCid_R(N)<\infty$.
\end{corollary}
\begin{proof}
We have that $R$ is Cohen-Macaulay since it admits a non-zero finitely generated $R$-module of finite injective dimension, so  $\Rfd_R (M) = \depth (R) - \depth (M)$ by \cite[Lemma 3.9]{FinitenessOfHomologicalDimensionOfExtModules}. By \cite[Corollary 3.7(a)]{CompleteIntersectionDimensionsandFoxbyClasses}, we may assume that $R$ is complete. Moreover, applying \cite[Corollary 3.17]{FinitenessOfHomologicalDimensionOfExtModules}, we have $\Gid_R (N)<\infty$. Therefore, we must then have that  $\Ext_R^i(M,N)=0$ for all $i>\depth (R)-\depth (M)$ by \cite[Theorem 2.10]{sazeedeh2013gorenstein} and $\RHom_R(M,N) \in \Dbf(R)$. Follows then by Lemma \ref{CI-homol}(1) that $\cdim_R(\RHom_R(M,N))< \infty$. Therefore, Theorem \ref{reqs22} gives the desired conclusion. 
\end{proof}

\noindent{\bf Acknowledgments.}  The authors thank Kaito Kimura for helpful discussions about his paper \cite{FinitenessOfHomologicalDimensionOfExtModules}, and in particular for sharing with us Lemmas \ref{triangle} and \ref{KimLemma}. The authors also wish to thank Lars Christensen for discussions about proving the isomorphisms in Theorems~\ref{issos} and ~\ref{HomGidInj}. The first author was supported by grants 2022/12114-0 and 2024/17809-1, S\~ao Paulo Research Foundation (FAPESP). The second author was supported by grants 2022/03372-5 and 2023/15733-5, S\~ao Paulo Research Foundation (FAPESP),  and by Dale M. Jensen Chair Foundation of the University of Nebraska–Lincoln. Part of this work was done when the second author was visiting the University of Nebraska-Lincoln;  he is grateful for the hospitality. Part of this work was done while the first author was visiting the University of Texas at Arlington. He is grateful for their kind hospitality. The first author also thanks the Department of Mathematics at the University of Nebraska-Lincoln for its hospitality during his visit in November 2025. 
\bibliographystyle{amsplain}
\bibliography{references.bib}
\end{document}